\documentclass[10pt,portrait,reqno,12pt]{amsart}%
\usepackage{amssymb,amsthm,amsmath, amsfonts, float}
\usepackage{graphicx}
\usepackage{multirow}
\usepackage{subfig}
\usepackage{amsmath}
\usepackage{amsfonts}
\usepackage{lscape}
\usepackage{amssymb}%
\providecommand{\U}[1]{\protect\rule{.1in}{.1in}}
\newtheorem{theorem}{Theorem}
\theoremstyle{plain}

\newtheorem{corollary}{Corollary}

\newtheorem{definition}{Definition}

\newtheorem{proposition}{Proposition}
\newtheorem{remark}{Remark}

\numberwithin{equation}{section}
\begin{document}
\title[The Gauss map of Rotational Hypersurfaces in $\mathbb E^4_1$]{Rotational Hypersufaces in $\mathbb E_1
^{4}$ with Generalized $L_{k}$ 1-Type Gauss Map}
\author[A. Kazan, M. Alt\i n and N.C. Turgay]{\bfseries Ahmet Kazan$^{1\ast}$, Mustafa Alt\i n$^{2}$ and Nurettin Cenk
Turgay$^{3}$}
\address{ \newline
$^{1}$\textit{Department of Computer Technologies, Do\u{g}an\c{s}ehir Vahap
K\"{u}\c{c}\"{u}k Vocational School, Malatya Turgut \"{O}zal University,
Malatya, Turkey} \newline
$^{2}$\textit{Technical Sciences Vocational School, Bing\"{o}l University,
Bing\"{o}l, Turkey} \newline
$^{3}$\textit{Department of Mathematics, Faculty of Science and Letters,
\.{I}stanbul Technical University, \.{I}stanbul, Turkey}\newline
$^{\ast}$\textit{Corresponding author: ahmet.kazan@ozal.edu.tr}}

\begin{abstract}
In this paper, we study the Gauss map of rotational hypersurfaces
in 4-dimensional Lorentz-Minkowski space concerning the linear second order differential operators
$L_{1}$  and $L_{2}$, where $L_1$ is usually called as the Cheng-Yau operator. We obtain  some classifications of rotational hypersurfaces which have
$L_{k}$-harmonic Gauss map, $L_{k}$-pointwise 1-type Gauss map and generalized $L_{k}%
$ 1-type Gauss map, where $k=1,2$.

\end{abstract}
\maketitle


\section{{Introduction}}

A crucial area of research in differential geometry is the study of
submanifolds of semi-Riemannian space forms whose coordinate functions are eigenfunctions of the Laplace
operator. Takahashi has
examined isometric immersions in the Euclidean space whose coordinate
functions are eigenfunctions of the Laplacian operator, associated to the same
eigenvalue $\lambda$ in \cite{Takahashi}, where it was proved that if an isometric immersion
$x:M\longrightarrow \mathbb R^{n+k}$ of a Riemannian $n$-manifold $M$ in a Euclidean
$(n+k)$-space satisfies $\Delta x=\lambda x$ for some constant $\neq0$, then
$\lambda$ must necessarily be positive and $x$ defines a minimal immersion in a
sphere $S^{n+k+1}$ of a radius $\sqrt{n/k}$ in $\mathbb E^{n+k}$. After this study, Garay has established an
extension of Takahashi's theorem in 1990 and he has proved that if
$x:M^{n}\longrightarrow \mathbb R^{n+1}$ is an immersed hypersurface in Euclidean
space with $\Delta x=Ax$, where $A\in \mathbb R^{(n+1)\times(n+1)}$ is a constant
diagonal matrix, then $M$ is minimal hypersurface in $\mathbb R^{n+1}$ or an open
piece of a round hypersphere or an open piece of generalized right spherical
cylinders \cite{Garay}. Also, Dillen et al. have considered surfaces in
$\mathbb R^{3}$ whose immersion satisfes the condition $\Delta x=Ax+b$, where $A\in
\mathbb R^{3\times3}$ is a constant matrix and $b\in \mathbb R^{3}$ is a constant vector in
\cite{Dillen}. Afterward, using different techniques, the result of Dillen et
al. has been independently generalized to hypersurfaces in $\mathbb R^{n+1}$ in
\cite{Chen} and \cite{Hasanis}.

The Laplace operator $\Delta$ of a hypersurface $M$ embedded in $\mathbb E^{n+1}$ is an
intrinsically defined second-order linear differential operator which
  can be viewed as the first operator in the
sequence of  operators $L_{0}=\Delta,\ L_1,\hdots,\ L_{n-1}$, where $L_{k}:C^\infty(M)\to C^\infty(M)$ represents the linearized operator of the first
variation of the $(k+1)$-th mean curvature, which arises from the normal
variations of the hypersurface defined by 
$$L_{k}(f)=tr(P_{k}\circ\nabla^{2}f),$$
where $\nabla^{2}f$  is the hessian of $f$ and  $P_{k}$ denotes the $k$-th Newton transformation associated with
the shape operator of the hypersurface.
Note that the operator $L_{1}$ is usually called as the Cheng-Yau operator, \cite{Alias2,Cheng,Reilly,Rosenberg}. In this context, an extension of
Takahashi theorem for the linearized operators $L_k$ has been given in \cite{Alias}, where  it is proved that a hypersurface of the Euclidean space $\mathbb E^{n+1}$ satisfying
$L_kx=Rx$ is an open part of a sphere $\mathbb S^{n}$. Later, this study was moved into the hypersurfaces of non-flat Riemannian space forms in in \cite{Yang2}.

On the other hand, the study of Gauss map of submanifolds by considering eigenfunctions of  the linearized operators $L_k$ have initiated by   Chen and Piccini in \cite{Chen-Piccinni}, where they consider submanifolds with 1-type Gauss map. After that, some important classes of hypersurfaces of semi-Euclidean spaces were studied in many papers by concerning the type of their Gauss maps. For example,  in \cite{Qian}
surfaces of revolution and canal surfaces with
generalized Cheng-Yau 1-type Gauss map have been obtained. Moreover, 
hypersurfaces in non-flat pseudo-Riemannian space forms satisfying  have
been studied in \cite{Lucas}. Also, pointwise 1-type Gauss map of surfaces in
$\mathbb E_{1}^{3}$ concerning Cheng-Yau operator has been given in \cite{Kim},
spacelike hypersurfaces in the Lorentz-Minkowski space satisfying
$L_kx=Rx+b$ have been classified in \cite{Yang}, hypersurfaces in non-flat
Lorentzian space forms satisfying $L_{k}\psi=A\psi+b$ have been studied in
\cite{Lucas2} and normalized null hypersurfaces in nonflat Lorentzian space
forms satisfying $L_kx=Ux+b$ have been introduced in \cite{Tetsing}.
Cheng-Yau operator and Gauss map of surfaces of revolution in $\mathbb E^{3}$ and
Cheng--Yau operator and Gauss map of rotational hypersurfaces in $\mathbb E^{4}$
have been studied in \cite{Kim2} and \cite{Guler}, respectively. Rotational
surfaces which have $L_{1}$-pointwise $1$-type Gauss map in Galilean 3-space
have been classified in \cite{Kelleci}.

In this paper, we give some classifications for Gauss maps of
rotational hypersurfaces in 4-dimensional Lorentz-Minkowski space concerning
the Cheng-Yau operator $L_{1}$ and as well as $L_{2}$. In Sect. 2,
 we firstly recall some fundamental notions about the linearized operators $L_{k}$
in $(n+1)$-dimensional Lorentz-Minkowski space and some geometric invariants
of rotational hypersurfaces in $4$-dimensional Lorentz-Minkowski space. In Sect. 3, we consider the Gauss map of rotational hypersurfaces by considering the Cheng-Yau operator. In Sect. 4, we obtain classifications of rotational hypersurfaces in $\mathbb E_{1}^{4}$ with generalized $L_2$ 1-type Gauss map, first and second kind $L_2$-pointwise 1-type Gauss map and $L_2$-harmonic Gauss map.

\section{{Preliminaries}}

Let $\mathbb E_{1}^{n+1}$ denote the $(n+1)$-dimensional Lorentz-Minkowski space with the
metric tensor given by 
$$\tilde g=\langle\,,\,\rangle=-dx_{1}^{2}+\sum\limits_{i=2}^ndx_{i}^{2},$$
where $(x_1,x_2,\hdots,x_n)$ is a Cartesian coordinate system in $\mathbb R^{n+1}$.

\subsection{Linearized Operators for the Gauss Map of Hypersurfaces in $\mathbb E_{1}^{n+1}$}

\

Consider an isometric immersion  $\Gamma:(M,g)\longrightarrow \mathbb E_{1}^{n+1}$ from an
$n$-dimensional orientable manifold $M$ into $\mathbb E_{1}^{n+1}$, where $g$ is a non-degenerate metric on $M$. Now, let $N$ denote the unit normal vector field and put $\left\langle N,N\right\rangle
=\varepsilon=\pm1$. Note that if  $\varepsilon=1$ (resp. $\varepsilon=-1$), then $g$ is Lorentzian (resp.  Riemannian). Moreover, in  \cite{ChengYauEn1} the equation
\begin{equation}
L_{k}N=-\varepsilon C_{k}\Big(  \nabla H_{k+1}+\left(  nH_{1}H_{k+1}-\left(
n-k-1\right)  H_{k+2}\right)  N\Big)   \label{Lk}%
\end{equation}
is obtained, where $H_{k}$ is the $k$-th mean curvature  of $M$,
defined by%
\begin{equation}
\left(
\begin{array}
[c]{c}%
n\\
k
\end{array}
\right)  H_{k}=\left(  -\varepsilon\right)  ^{k}a_{k} , \label{Hk}%
\end{equation}
and the constants $a_{k}$ and $C_{k}$ are given by
\begin{equation}
\left.
\begin{array}
[c]{l}%
a_{1}=-\sum_{i=1}^{n}\kappa_{i},\\
\\
a_{k}=(-1)^{k}\sum_{i_{1}<i_{2}<...<i_{k}}^{n}\kappa_{i_{1}}\kappa_{i_{2}%
}...\kappa_{i_{k}},~k=2,3,...,n
\end{array}
\right\}  \label{ak}%
\end{equation}
and 
\begin{equation}
C_{k}=\left(
\begin{array}
[c]{c}%
n\\
k+1
\end{array}
\right)  (-\varepsilon)^{k}. \label{Ck}%
\end{equation}

\begin{definition}
Let $\mathfrak{m}$ and $\mathfrak{n}$ be non-zero smooth functions on $M$, $C\in
\mathbb E_{1}^{n+1}$ be a non-zero constant vector and $k\in\{0,1,2,...,n-1\}$. If the
Gauss map $N$ of an oriented submanifold $M$ in $\mathbb E_{1}^{4}$ satisfies

\begin{description}
\item[i] $L_{k}N=0$, then $N$ is called $L_{k}$-harmonic;

\item[ii] $L_{k}N=\mathfrak{m}N,$ then $M$ has first kind $L_{k}$-pointwise
1-type Gauss map;

\item[iii] $L_{k}N=\mathfrak{m}(N+C),$ then $M$ has second kind $L_{k}%
$-pointwise 1-type Gauss map;

\item[iv] $L_{k}N=\mathfrak{m}N+\mathfrak{n}C,$ then $M$ has generalized
$L_{k}$ 1-type Gauss map.
\end{description}
\end{definition}

Now, consider the case $n=3$ and let  $M$ be an oriented hypersurface in $\mathbb E_{1}^{4}$  with  a local coordinat system $\{s,t,w\}$. Then, the gradient of a smooth function $f\in C^\infty (M)$ is defined by 
\begin{equation}
\nabla f=\frac{1}{\mathfrak{g}}
\left(  \nabla_1  \partial_s+ \nabla_2\partial_t+\nabla_3 \partial_w
\right)  , \label{gradf}%
\end{equation}
where we put 
\begin{eqnarray*}
\mathfrak{g}&=&g_{13}^{2}g_{22}-2g_{12}g_{13}g_{23}+g_{11}g_{23}^{2}+g_{12}%
^{2}g_{33}-g_{11}g_{22}g_{33},\\
\nabla_1&=& \left( g_{23}^{2}-g_{22}g_{33}\right)  f_{s}+\left(  -g_{13}%
g_{23}+g_{12}g_{33}\right)  f_{t}+\left(  g_{13}g_{22}-g_{12}g_{23}\right)
f_{w},\\
\nabla_2&=&    \left(  -g_{13}g_{23}+g_{12}g_{33}\right)  f_{s}+\left(  g_{13}%
^{2}-g_{11}g_{33}\right)  f_{t}+\left(  -g_{12}g_{13}+g_{11}g_{23}\right)
f_{w} ,\\
\nabla_3&=&  \left(  g_{13}g_{22}-g_{12}g_{23}\right)  f_{s}+\left(  -g_{12}%
g_{13}+g_{11}g_{23}\right)  f_{t}+\left(  g_{12}^{2}-g_{11}g_{22}\right)
f_{w}  
\end{eqnarray*}
and $g_{11}=\left\langle \partial s,\partial
s\right\rangle ,$ $g_{12}=\left\langle \partial s,\partial t\right\rangle ,$
$g_{13}=\left\langle \partial s,\partial w\right\rangle ,$ $g_{22}%
=\left\langle \partial t,\partial t\right\rangle ,$ $g_{23}=\left\langle
\partial t,\partial w\right\rangle ,$ $g_{33}=\left\langle \partial w,\partial
w\right\rangle $.

\subsection{Rotational hypersurfaces in $\mathbb E_{1}^{4}$}

\

In this subsection we recall some geometric characterizations of rotational
hypersurfaces in $\mathbb E_{1}^{4}$. Note that there are three different classes of rotational hypersurfaces in $\mathbb E_{1}^{4}$ subject to the causality of their axis of rotation (for more datails, see \cite{hacettepe}).

The rotational hypersurface $\Gamma^{s}$ which is obtained by rotating the
profile curve $\alpha^{s}(s)=(s,0,0,f(s))$ about spacelike axis $(0,0,0,1)$ is
given by%
\begin{equation}
\Gamma^{s}(s,t,w)=\left(  s\cosh t\cosh w,s\sinh t,s\cosh t\sinh
w,f(s)\right)  , \label{surfspace}%
\end{equation}
where $s\in%
\mathbb{R}
-\{0\}$.

With the aid of the first differentials of (\ref{surfspace}) with respect to
$s,$ $t$ and $w$, the Gauss map $N^{s}$ of the rotational hypersurface
(\ref{surfspace}) is obtained by
\begin{equation}
N^{s}=-\frac{1}{\sqrt{1-f^{\prime2}}}\left(  f^{\prime}\cosh t\cosh
w,f^{\prime}\sinh t,f^{\prime}\cosh t\sinh w,1\right)  \label{normalspace}%
\end{equation}
and so,%
\begin{equation}
\varepsilon^{s}=\left\langle N^{s},N^{s}\right\rangle =1. \label{epsilonspace}%
\end{equation}
Also, the nonzero components of the coefficients in the first fundamental form
of (\ref{surfspace}) are%
\begin{equation}
g_{11}^{s}=f^{\prime2}-1,\text{ }g_{22}^{s}=s^{2},\text{ }g_{33}^{s}%
=s^{2}\cosh^{2}t \label{gspace}%
\end{equation}
and the principal curvatures of (\ref{surfspace}) are%

\begin{equation}
\kappa_{1}^{s}=\frac{f^{\prime\prime}}{\left(  1-f^{\prime2}\right)
^{\frac{3}{2}}},\text{ }\kappa_{2}^{s}=\kappa_{3}^{s}=\frac{f^{\prime}}%
{s\sqrt{1-f^{\prime2}}}. \label{aslispace}%
\end{equation}

On the other hand, the rotational hypersurface $\Gamma^{t}$ which is obtained by rotating the
profile curve $\alpha^{t}(s)=(f(s),0,0,s)$ about timelike axis $(1,0,0,0)$ is
given by%
\begin{equation}
\Gamma^{t}(s,t,w)=\left(  f(s),-s\cos t\sin w,-s\sin t,s\cos t\cos w\right)  ,
\label{surftime}%
\end{equation}
where $s\in%
\mathbb{R}
-\{0\},$ $0\leq t,w\leq2\pi$.

Moreover, the Gauss map $N^{t}$ of the rotational hypersurface
(\ref{surftime}) is defined by
\begin{equation}
N^{t}=\frac{1}{\sqrt{f^{\prime2}-1}}\left(  1,-f^{\prime}\cos t\sin
w,-f^{\prime}\sin t,f^{\prime}\cos t\cos w\right)  \label{normaltime}%
\end{equation}
and so%
\begin{equation}
\varepsilon^{t}=\left\langle N^{t},N^{t}\right\rangle =1. \label{epsilontime}%
\end{equation}
Also, the nonzero components of the coefficients in the first fundamental form
of (\ref{surftime}) are%
\begin{equation}
g_{11}^{t}=1-f^{\prime2},\text{ }g_{22}^{t}=s^{2},\text{ }g_{33}^{t}=s^{2}%
\cos^{2}t \label{gtime}%
\end{equation}
and the principal curvatures of (\ref{surftime}) are%

\begin{equation}
\kappa_{1}^{t}=\frac{f^{\prime\prime}}{\left(  f^{\prime2}-1\right)
^{\frac{3}{2}}},\text{ }\kappa_{2}^{t}=\kappa_{3}^{t}=-\frac{f^{\prime}%
}{s\sqrt{f^{\prime2}-1}}. \label{aslitime}%
\end{equation}

Finally, the rotational hypersurface $\Gamma^{l}$ which is obtained by rotating the
profile curve $\alpha^{l}(s)=(s,f(s),0,0)$ about lightlike axis $(1,1,0,0)$ is
given by%
\begin{equation}
	\Gamma^{l}(s,t,w)=%
	\begin{array}
		[c]{l}%
		\left(  \left(  \frac{t^{2}+w^{2}}{2}+1\right)  s-\frac{t^{2}+w^{2}}%
		{2}f(s),\frac{t^{2}+w^{2}}{2}s+\left(  1-\frac{t^{2}+w^{2}}{2}\right)
		f(s),\right.  \\
		\text{\ \ \ \ \ \ \ \ \ \ \ \ \ \ \ \ \ \ \ \ \ \ \ \ \ \ \ \ \ \ \ \ \ \ \ \ \ \ \ \ \ }st-f(s)t,sw-f(s)w \big),
	\end{array}
	\label{surflight}%
\end{equation}
where $s\in%
\mathbb{R}
-\{0\}$.

It turns out that the Gauss map $N^{l}$ of the rotational hypersurface
(\ref{surflight}) is 
\begin{equation}
	N^{l}=\frac{1}{2\sqrt{1-f^{\prime2}}}%
	\begin{array}
		[c]{l}%
		\big(  t^{2}+w^{2}-\left(  t^{2}+w^{2}+2\right)  f^{\prime},t^{2}%
		+w^{2}-2-\left(  t^{2}+w^{2}\right)  f^{\prime},  \\
		  \text{\ \ \ \ \ \ \ \ \ \ \ \ \ \ \ \ \ \ \ \ \ \ \ \ \ \ \ \ \ \ \ \ \ \ \ \ \ \ \ }2t\left(  1-f^{\prime}\right)  ,2w\left(  1-f^{\prime}\right)
		\big)
	\end{array}
	\label{normallight}%
\end{equation}
and so%
\begin{equation}
\varepsilon^{l}=\left\langle N^{l},N^{l}\right\rangle =1. \label{epsilonlight}%
\end{equation}
Also, the nonzero components of the coefficients in the first fundamental form
of (\ref{surflight}) are%
\begin{equation}
g_{11}^{l}=f^{\prime2}-1,\text{ }g_{22}^{l}=\left(  s-f\right)  ^{2},\text{
}g_{33}^{l}=\left(  s-f\right)  ^{2} \label{glight}%
\end{equation}
and the principal curvatures of (\ref{surflight}) are%

\begin{equation}
\kappa_{1}^{l}=\frac{f^{\prime\prime}}{\left(  1-f^{\prime2}\right)
^{\frac{3}{2}}},\text{ }\kappa_{2}^{l}=\kappa_{3}^{l}=\frac{f^{\prime}%
-1}{\left(  s-f\right)  \sqrt{1-f^{\prime2}}}. \label{aslilight}%
\end{equation}

\begin{remark} 
Throughout this study, we are going to put $f=f(u),$  and $\prime$ is going to stand for the ordinary derivative with respect to $u$, i.e., 
$$f^{\prime}%
=\frac{df }{du}.$$ 
Also, $c_{i}$ will denote some real constants for $i\in \mathbb{N^+}$.
\end{remark}
\section{On Hypersurfaces with Generalized $L_1$ 1-Type Gauss Map}

In this section, we obtain the $L_{1}$ (Cheng-Yau) operator of the Gauss maps
of the rotational hypersurfaces (\ref{surfspace}), (\ref{surftime}),
(\ref{surflight}) and give some classifications for these hypersurfaces which have $L_{1}$-harmonic Gauss map, first kind $L_{1}$-pointwise 1-type Gauss map, second kind $L_{1}$-pointwise 1-type Gauss map and generalized $L_{1}$ 1-type Gauss map, seperately.

\subsection{Rotational Hypersurfaces about Spacelike Axis}

\

The functions $a_{k}$ of the rotational hypersurface (\ref{surfspace}) in
$\mathbb E_{1}^{4}$ are obtained from (\ref{ak}) and (\ref{aslispace}) by%
\begin{equation}
\left.
\begin{array}
[c]{l}%
a_{1}^{s}=\frac{2f^{\prime}\left(  f^{\prime2}-1\right)  -sf^{\prime\prime}%
}{s\left(  1-f^{\prime2}\right)  ^{3/2}},\\
\\
a_{2}^{s}=\frac{f^{\prime}\left(  2sf^{\prime\prime}+f^{\prime}\left(
1-f^{\prime2}\right)  \right)  }{s^{2}\left(  1-f^{\prime2}\right)  ^{2}},\\
\\
a_{3}^{s}=\frac{-f^{\prime2}f^{\prime\prime}}{s^{2}\left(  1-f^{\prime
2}\right)  ^{5/2}}.
\end{array}
\right\}  \label{akspace}%
\end{equation}
Also, from (\ref{gradf}), (\ref{gspace}) and (\ref{akspace}), we have%
\begin{equation}
\nabla a_{2}^{s}=\frac{2P_{1}^{s}}{s^{3}\left(  f^{\prime2}-1\right)  ^{4}%
}\left(  \cosh t\cosh w,\sinh t,\cosh t\sinh w,f^{\prime}\right)  ,
\label{a2nablaspace}%
\end{equation}
where
\[
P_{1}^{s}=f^{\prime}\left(  s^{2}\left(  f^{\prime2}-1\right)  f^{\prime
\prime\prime}+\left(  \left(  f^{\prime2}-1\right)  ^{2}-3s^{2}f^{\prime
\prime2}\right)  \right)  f^{\prime}-s^{2}f^{\prime\prime2}.
\]

Therefore, (\ref{Lk}), (\ref{Hk}), (\ref{normalspace}), (\ref{epsilonspace}),
(\ref{akspace}) and (\ref{a2nablaspace}) imply
\begin{equation}
L_{1}N^{s}=\frac{-2}{s^{3}\left(  f^{\prime2}-1\right)  ^{4}}\left(  \cosh
t\cosh wQ_{1}^{s},\sinh tQ_{1}^{s},\cosh t\sinh wQ_{1}^{s},R_{1}^{s}\right)
, \label{L1space}%
\end{equation}
where
\begin{align*}
Q_{1}^{s}  &  =3\left(  1-f^{\prime2}\right)  f^{\prime4}+f^{\prime
8}-sf^{\prime5}f^{\prime\prime}+s^{2}f^{\prime\prime2}+\left(  4s^{2}%
f^{\prime\prime2}-1\right)  f^{\prime2}+s^{2}f^{\prime}f^{\prime\prime\prime
}+s\left(  f^{\prime\prime}-sf^{\prime\prime\prime}\right)  f^{\prime3},\\
R_{1}^{s}  &  =-sf^{\prime}\left(  \left(  \left(  f^{\prime2}-1\right)
f^{\prime}-s\left(  3f^{\prime2}+2\right)  f^{\prime\prime}\right)
f^{\prime\prime}+s\left(  f^{\prime2}-1\right)  f^{\prime}f^{\prime
\prime\prime}\right)  .
\end{align*}

If the rotational hypersurface (\ref{surfspace}) in $\mathbb E_{1}^{4}$ is flat (resp.
minimal), then $a_{3}^{s}=0$ (resp. $a_{1}^{s}=0$)$.$ So, from (\ref{L1space}%
), we get

\begin{proposition}
If the rotational hypersurface (\ref{surfspace}) in $\mathbb E_{1}^{4}$ is flat (resp.
minimal), then we obtain%
\begin{align}
&  L_{1}N^{s}=\frac{-1}{s^{3}\left(  f^{\prime2}-1\right)  ^{4}}\left(  \cosh
t\cosh wA_{1}^{s},\sinh tA_{1}^{s},\cosh t\sinh wA_{1}^{s},A_{2}^{s}\right)
\label{L1spaceflat}\\
&  \left(  \text{resp. }L_{1}N^{s}=\frac{1}{s^{3}\left(  f^{\prime2}-1\right)
^{4}}\left(  \cosh t\cosh wA_{3}^{s},\sinh tA_{3}^{s},\cosh t\sinh wA_{3}%
^{s},A_{4}^{s}\right)  \right)  , \label{L1spaceminimal}%
\end{align}
where
\begin{align*}
&  A_{1}^{s}=6\left(  1-f^{\prime2}\right)  f^{\prime4}+2f^{\prime
8}-5sf^{\prime5}f^{\prime\prime}+2s^{2}f^{\prime\prime2}+\left(
8s^{2}f^{\prime\prime2}-2\right)  f^{\prime2}+2s^{2}f^{\prime}f^{\prime
\prime\prime}+s\left(  5f^{\prime\prime}-2sf^{\prime\prime\prime}\right)
f^{\prime3},\\
&  A_{2}^{s}=-s\left(  2s\left(  f^{\prime2}-1\right)  f^{\prime}%
f^{\prime\prime\prime}+\left(  5\left(  f^{\prime2}-1\right)  f^{\prime
}-2s\left(  3f^{\prime2}+2\right)  f^{\prime\prime}\right)  f^{\prime\prime
}\right)  f^{\prime},\\
&  \left(
\begin{array}
[c]{c}%
\text{resp. }A_{3}^{s}=-3s\left(  f^{\prime2}-1\right)  f^{\prime3}%
f^{\prime\prime}-2s^{2}f^{\prime\prime2}+2\left(  \left(  \left(  f^{\prime
2}-1\right)  ^{2}-3s^{2}f^{\prime\prime2}\right)  f^{\prime}+s^{2}\left(
f^{\prime2}-1\right)  f^{\prime\prime\prime}\right)  f^{\prime},\\
A_{4}^{s}=\left(  -2\left(  2-f^{\prime2}\right)  f^{\prime4}-2s^{2}%
f^{\prime\prime2}+\left(  2-6s^{2}f^{\prime\prime2}\right)  f^{\prime
2}+s\left(  3f^{\prime\prime}-2sf^{\prime\prime\prime}\right)  f^{\prime
}+s\left(  2sf^{\prime\prime\prime}-3f^{\prime\prime}\right)  f^{\prime
3}\right)  f^{\prime}.
\end{array}
\right)
\end{align*}

\end{proposition}

\begin{remark}
\label{remarkspace}We know that \cite{hacettepe}, if the rotational
hypersurface (\ref{surfspace}) in $\mathbb E_{1}^{4}$ is flat (resp. minimal), then we
have $f(s)=\frac{c_{1}{}^{1/3}s}{\sqrt{1+c_{1}{}^{2/3}}}$ (resp. $f^{\prime
}(s)=\frac{c_{2}}{\sqrt{s^{4}+c_{2}{}^{2}}}$).
\end{remark}

Thus, from (\ref{L1space}) or (\ref{L1spaceflat}) (resp. (\ref{L1spaceminimal}%
)),\ then we have

\begin{corollary}
Let $M$ be the rotational hypersurface given by (\ref{surfspace}) in $\mathbb E_{1}^{4}$. Then, we have the followings:
\begin{enumerate}
\item[(i)] If $M$ is flat, then
\begin{equation}\nonumber
  L_{1}N^{s}=\frac{2c_{1}{}^{2/3}}{s^{3}}\left(  \cosh t\cosh w,\sinh
t,\cosh t\sinh w,0\right).
\end{equation}
\item[(ii)] If $M$ is minimal, then
\begin{align*}
\begin{split}
L_{1}N^{s}=\frac{-6c_{2}{}^{2}}{s^{11}}&\left(
\left(  3s^{4}+4c_{2}{}^{2}\right)  \cosh t\cosh w,\left(  3s^{4}+4c_{2}{}%
^{2}\right)  \sinh t,\right.\\ &\left.
\left(  3s^{4}+4c_{2}{}^{2}\right)  \cosh t\sinh w,4c_{2}\sqrt{s^{4}+c_{2}%
{}^{2}}
\right). 
\end{split}
\end{align*}
\end{enumerate}
\end{corollary}

Now, let us give some classifications for the rotational hypersurface
(\ref{surfspace}) which have $L_{1}$-harmonic Gauss map, first kind $L_{1}%
$-pointwise 1-type Gauss map, second kind $L_{1}$-pointwise 1-type Gauss map
and generalized $L_{1}$ 1-type Gauss map.

If $\Gamma^{s}(s,t,w)$ has a generalized $L_{1}$ (Cheng-Yau) 1-type Gauss map,
i.e., $L_{1}N^{s}=\mathfrak{m}N^{s}+\mathfrak{n}C,$ where $C=\left(
C_{1},C_{2},C_{3},C_{4}\right)  $ is a constant vector, then from
(\ref{normalspace}) and (\ref{L1space}), we get%
\begin{equation}
\left.
\begin{array}
[c]{l}%
\frac{-2\cosh t\cosh wQ_{1}^{s}}{s^{3}\left(  f^{\prime2}-1\right)  ^{4}%
}=\mathfrak{m}\left(  \frac{-f^{\prime}\cosh t\cosh w}{\sqrt{1-f^{\prime2}}%
}\right)  +\mathfrak{n}C_{1},\\
\\
\frac{-2\sinh tQ_{1}^{s}}{s^{3}\left(  f^{\prime2}-1\right)  ^{4}%
}=\mathfrak{m}\left(  \frac{-f^{\prime}\sinh t}{\sqrt{1-f^{\prime2}}}\right)
+\mathfrak{n}C_{2},\\
\\
\frac{-2\cosh t\sinh wQ_{1}^{s}}{s^{3}\left(  f^{\prime2}-1\right)  ^{4}%
}=\mathfrak{m}\left(  \frac{-f^{\prime}\cosh t\sinh w}{\sqrt{1-f^{\prime2}}%
}\right)  +\mathfrak{n}C_{3},\\
\\
\frac{-2R_{1}^{s}}{s^{3}\left(  f^{\prime2}-1\right)  ^{4}}=\mathfrak{m}%
\left(  \frac{-1}{\sqrt{1-f^{\prime2}}}\right)  +\mathfrak{n}C_{4}.
\end{array}
\right\}  \label{sinif1space}%
\end{equation}
It is obvious from the first three equations of (\ref{sinif1space}), we get
$C_{1}=C_{2}=C_{3}=0$. Here if $f^{\prime}=0,$ then we have $L_{1}N^{s}=0.$ Thus,

\begin{theorem}
The rotational hypersurface
\[
\Gamma^{s}(s,t,w)=\left(  s\cosh t\cosh w,s\sinh t,s\cosh t\sinh
w,c_{3}\right)
\]
has $L_{1}$-harmonic Gauss map in $\mathbb E_{1}^{4}.$
\end{theorem}

From now on, we assume that $f^{\prime}\neq0.$ In (\ref{sinif1space}), we can
take $C_{4}=0$ or $C_{4}\neq0.$

Firstly let us suppose that $C_{4}=0.$ In this case; from the first three
equations and the last equation of (\ref{sinif1space}), we have%
\begin{equation}
\mathfrak{m}=-\frac{-2\left(  f^{\prime2}-1\right)  ^{3}f^{\prime2}+2s\left(
f^{\prime2}-1\right)  f^{\prime3}f^{\prime\prime}-2s^{2}\left(  4f^{\prime
2}+1\right)  f^{\prime\prime2}+2s^{2}\left(  f^{\prime2}-1\right)  f^{\prime
}f^{\prime\prime\prime}}{s^{3}f^{\prime}\left(  1-f^{\prime2}\right)  ^{7/2}}
\label{sinif2space}%
\end{equation}
and%
\begin{equation}
\mathfrak{m}=-\frac{2f^{\prime}\left(  s\left(  f^{\prime2}-1\right)
f^{\prime}f^{\prime\prime\prime}+\left(  \left(  f^{\prime2}-1\right)
f^{\prime}-s\left(  3f^{\prime2}+2\right)  f^{\prime\prime}\right)
f^{\prime\prime}\right)  }{s^{2}\left(  1-f^{\prime2}\right)  ^{7/2}},
\label{sinif3space}%
\end{equation}
respectively. From (\ref{sinif2space}) and (\ref{sinif3space}), we get%
\begin{equation}
f^{\prime}\left(  s^{2}\left(  f^{\prime2}-1\right)  f^{\prime\prime\prime
}+\left(  \left(  f^{\prime2}-1\right)  ^{2}-3s^{2}f^{\prime\prime2}\right)
f^{\prime}\right)  -s^{2}f^{\prime\prime2}=0. \label{sinif4space}%
\end{equation}
By solving (\ref{sinif4space}), we obtain that%
\[
f=\pm\int_{1}^{s}\frac{\sqrt{2\left(  6c_{4}c_{5}x^{3}+c_{4}\right)  }}%
{\sqrt{12c_{4}c_{5}x^{3}+3x+2c_{4}}}dx+c_{6}.
\]
Hence,

\begin{theorem}
The rotational hypersurface
\[
\Gamma^{s}(s,t,w)=\left(  s\cosh t\cosh w,s\sinh t,s\cosh t\sinh w,\pm\int%
_{1}^{s}\frac{\sqrt{2\left(  6c_{4}c_{5}x^{3}+c_{4}\right)  }}{\sqrt
{12c_{4}c_{5}x^{3}+3x+2c_{4}}}dx+c_{6}\right)
\]
has first kind $L_{1}$-pointwise 1-type Gauss map, i.e., $L_{1}N^{s}%
=\mathfrak{m}N^{s},$ in $\mathbb E_{1}^{4},$ where%
\begin{align*}
&  \mathfrak{m}=-\frac{c_{4}{}^{2}\sqrt{2}\left(  1+12c_{5}s^{3}\left(
1+12c_{5}s^{3}\right)  \right)  }{s^{9/2}\sqrt{3c_{4}\left(  1+6c_{5}%
s^{3}\right)  }}\\
&  \text{and}\\
&  L_{1}N^{s}=\frac{2\left(  1+12c_{5}s^{3}\left(  1+12c_{5}s^{3}\right)
\right)  }{-3c_{4}{}^{-2}s^{5}}\left(  \cosh t\cosh w,\sinh t,\cosh t\sinh
w,\frac{\sqrt{12c_{4}c_{5}s^{3}+3s+2c_{4}}}{-\sqrt{2c_{4}\left(  1+6c_{5}%
s^{3}\right)  }}\right)  .
\end{align*}

\end{theorem}

Now let us suppose that $C_{4}\neq0.$ In this case; from the first three
equations and the last equation of (\ref{sinif1space}), we have%
\begin{equation}
\left.
\begin{array}
[c]{l}%
\mathfrak{m}=-\frac{-2\left(  f^{\prime2}-1\right)  ^{3}f^{\prime2}+2s\left(
f^{\prime2}-1\right)  f^{\prime3}f^{\prime\prime}-2s^{2}\left(  \left(
4f^{\prime2}+1\right)  f^{\prime\prime2}-\left(  f^{\prime2}-1\right)
f^{\prime}f^{\prime\prime\prime}\right)  }{s^{3}f^{\prime}\left(
1-f^{\prime2}\right)  ^{7/2}}\\
\text{and}\\
\mathfrak{n}=\frac{2\left(  s^{2}f^{\prime\prime\prime}\left(  f^{\prime
2}-1\right)  +f^{\prime}\left(  \left(  f^{\prime2}-1\right)  ^{2}%
-3s^{2}f^{\prime\prime2}\right)  \right)  f^{\prime}-2s^{2}f^{\prime\prime2}%
}{C_{4}s^{3}f^{\prime}\left(  f^{\prime2}-1\right)  ^{3}},
\end{array}
\right\}  \label{sinif5space}%
\end{equation}
respectively.

\begin{theorem}
The rotational hypersurface (\ref{surfspace}) has generalized $L_{1}%
$ 1-type Gauss map, i.e., $L_{1}N^{s}=\mathfrak{m}N^{s}%
+\mathfrak{n}C,$ in $\mathbb E_{1}^{4},$ where $\mathfrak{m}$ and $\mathfrak{n}$ are
non-zero smooth functions given by (\ref{sinif5space}) and $C=(0,0,0,C_{4})$
is a non-zero vector with non-zero constant $C_{4}$.
\end{theorem}

If we assume $\mathfrak{m}=\mathfrak{n}$ in (\ref{sinif1space}) (or in
(\ref{sinif5space})), then we obtain%
\begin{align}\label{sinif7space}%
\begin{split}
0=&  \left(  C_{4}\sqrt{1-f^{\prime2}}-1\right)  \left(
2\left(  f^{\prime2}-1\right)  ^{3}f^{\prime2}-2s\left(  f^{\prime2}-1\right)
f^{\prime3}f^{\prime\prime}\right.\\&\left.
+2s^{2}\left(  4f^{\prime2}+1\right)  f^{\prime\prime2}-2s^{2}\left(
f^{\prime2}-1\right)  f^{\prime}f^{\prime\prime\prime}\right) \\
&  -2s\left(  s\left(  f^{\prime2}-1\right)  f^{\prime}f^{\prime\prime\prime
}+\left(  \left(  f^{\prime2}-1\right)  f^{\prime}-s\left(  3f^{\prime
2}+2\right)  f^{\prime\prime}\right)  f^{\prime\prime}\right)  f^{\prime2}.
\end{split}
\end{align}
Therefore,

\begin{theorem}
The rotational hypersurface (\ref{surfspace}) has second kind $L_{1}%
$-pointwise 1-type Gauss map, i.e., $L_{1}N^{s}=\mathfrak{m}\left(
N^{s}+C\right)  ,$ in $\mathbb E_{1}^{4},$ if the differential equation
(\ref{sinif7space}) holds.
\end{theorem}

\subsection{Rotational Hypersurfaces about Timelike Axis}

\

The functions $a_{k}$ of rotational hypersurface (\ref{surftime}) in
$\mathbb E_{1}^{4}$ are obtained from (\ref{ak}) and (\ref{aslitime}) by%
\begin{equation}
\left.
\begin{array}
[c]{l}%
a_{1}^{t}=\frac{2f^{\prime}\left(  f^{\prime2}-1\right)  -sf^{\prime\prime}%
}{s\left(  f^{\prime2}-1\right)  ^{3/2}},\\
\\
a_{2}^{t}=\frac{f^{\prime}\left(  f^{\prime}\left(  f^{\prime2}-1\right)
-2sf^{\prime\prime}\right)  }{s^{2}\left(  f^{\prime2}-1\right)  ^{2}},\\
\\
a_{3}^{t}=\frac{-f^{\prime2}f^{\prime\prime}}{s^{2}\left(  f^{\prime
2}-1\right)  ^{5/2}}.
\end{array}
\right\}  \label{aktime}%
\end{equation}
Also, from (\ref{gradf}), (\ref{gtime}) and (\ref{aktime}), we have%
\begin{equation}
\nabla a_{2}^{t}=\frac{2P_{1}^{t}}{s^{3}\left(  f^{\prime2}-1\right)  ^{4}%
}\left(  f^{\prime},-\cos t\sin w,-\sin t,\cos t\cos w\right),
\label{a2nablatime}%
\end{equation}
where
\[
P_{1}^{t}=-s^{2}f^{\prime\prime2}+\left(  \left(  \left(  f^{\prime
2}-1\right)  ^{2}-3s^{2}f^{\prime\prime2}\right)  f^{\prime}+s^{2}\left(
f^{\prime2}-1\right)  f^{\prime\prime\prime}\right)  f^{\prime}.
\]

So, from (\ref{Lk}), (\ref{Hk}), (\ref{normaltime}), (\ref{epsilontime}),
(\ref{aktime}) and (\ref{a2nablatime}), we have%
\begin{equation}
L_{1}N^{t}=\frac{2}{s^{3}\left(  f^{\prime2}-1\right)  ^{4}}\left(  R_{1}%
^{t},\cos t\sin wQ_{1}^{t},\sin tQ_{1}^{t},-\cos t\cos wQ_{1}^{t}\right)  ,
\label{L1time}%
\end{equation}
where
\begin{align*}
Q_{1}^{t}  &  =3\left(  1-f^{\prime2}\right)  f^{\prime4}+f^{\prime
8}-sf^{\prime5}f^{\prime\prime}+s^{2}f^{\prime\prime2}+\left(  4s^{2}%
f^{\prime\prime2}-1\right)  f^{\prime2}+s^{2}f^{\prime}f^{\prime\prime\prime
}+sf^{\prime3}\left(  f^{\prime\prime}-sf^{\prime\prime\prime}\right)  ,\\
R_{1}^{t}  &  =sf^{\prime}\left(  \left(  f^{\prime}\left(  f^{\prime
2}-1\right)  -s\left(  3f^{\prime2}+2\right)  \right)  f^{\prime\prime
}+s\left(  f^{\prime2}-1\right)  f^{\prime}f^{\prime\prime\prime}\right)  .
\end{align*}

Thus,

\begin{proposition}
If the rotational hypersurface (\ref{surftime}) in $\mathbb E_{1}^{4}$ is flat (resp.
minimal), then we obtain%
\begin{align}
&  L_{1}N^{t}=\frac{1}{s^{3}\left(  f^{\prime2}-1\right)  ^{4}}\left(
A_{2}^{t},\cos t\sin wA_{1}^{t},\sin tA_{1}^{t},-\cos t\cos wA_{1}^{t}\right)
\label{L1timeflat}\\
&  \left(  \text{resp. }L_{1}N^{t}=\frac{1}{s^{3}\left(  f^{\prime2}-1\right)
^{4}}\left(  A_{4}^{t},\cos t\sin wA_{3}^{t},\sin tA_{3}^{t},-\cos t\cos
wA_{3}^{t}\right)  ,\right)  \label{L1timeminimal}%
\end{align}
where
\begin{align*}
&  A_{1}^{t}=6f^{\prime4}\left(  1-f^{\prime2}\right)  +2f^{\prime
8}-5sf^{\prime5}f^{\prime\prime}+2s^{2}f^{\prime\prime2}+\left(
8s^{2}f^{\prime\prime2}-2\right)  f^{\prime2}+2s^{2}f^{\prime}f^{\prime
\prime\prime}+sf^{\prime3}\left(  5f^{\prime\prime}-2sf^{\prime\prime\prime
}\right)  ,\\
&  A_{2}^{t}=sf^{\prime}\left(  \left(  5f^{\prime}\left(  f^{\prime
2}-1\right)  -2s\left(  3f^{\prime2}+2\right)  f^{\prime\prime}\right)
f^{\prime\prime}+2s\left(  f^{\prime2}-1\right)  f^{\prime}f^{\prime
\prime\prime}\right)  ,\\
&  \left(
\begin{array}
[c]{c}%
\text{resp. }A_{3}^{t}=3s\left(  f^{\prime2}-1\right)  f^{\prime3}%
f^{\prime\prime}+2s^{2}f^{\prime\prime2}-2\left(  f^{\prime}\left(  \left(
f^{\prime2}-1\right)  ^{2}-3s^{2}f^{\prime\prime2}\right)  +s^{2}\left(
f^{\prime2}-1\right)  f^{\prime\prime\prime}\right)  f^{\prime},\\
A_{4}^{t}=f^{\prime}\left(  -2f^{\prime4}\left(  2-f^{\prime2}\right)
-2s^{2}f^{\prime\prime2}+\left(  2-6s^{2}f^{\prime\prime2}\right)  f^{\prime
2}+\left(  3f^{\prime\prime}-2sf^{\prime\prime\prime}\right)  sf^{\prime
}\left(  1-f^{\prime2}\right)  \right)  .
\end{array}
\right)
\end{align*}

\end{proposition}

\begin{remark}
\label{remarktime}We know that \cite{hacettepe}, if the rotational
hypersurface (\ref{surftime}) in $\mathbb E_{1}^{4}$ is flat (resp. minimal), then we
have $f(s)=\frac{c_{7}{}^{1/3}s}{\sqrt{c_{7}{}^{2/3}-1}}$ (resp. $f^{\prime
}(s)=\frac{c_{8}}{\sqrt{-s{}^{4}+c_{8}{}^{2}}}$).
\end{remark}

Thus, from (\ref{L1time}) or (\ref{L1timeflat}) (resp. (\ref{L1timeminimal})),
we have

\begin{corollary}
Let $M$ be the rotational hypersurface given by (\ref{surftime}) in $\mathbb E_{1}^{4}$. Then, we have the followings: 
\begin{enumerate}
\item[(i)] If $M$ is flat, then
\begin{equation}\nonumber
L_{1}N^{t}=\frac{2c_{7}{}^{2/3}}{s^{3}}\left(  0,\cos t\sin w,\sin t,-\cos
t\cos w\right).
\end{equation}
\item[(ii)] If $M$ is minimal, then
\begin{align*}
\begin{split}
L_{1}N^{t}=\frac{6c_{8}{}^{2}}{s^{11}}&\left(
-4c_{8}{}\sqrt{-s^{4}+c_{8}{}^{2}},\left(  -3s^{4}+4c_{8}{}^{2}\right)  \cos
t\sin w,\right.\\&\left.
\left(  -3s^{4}+4c_{8}{}^{2}\right)  \sin t,\left(  3s^{4}-4c_{8}{}%
^{2}\right)  \cos t\cos w
\right). 
\end{split}
\end{align*}
\end{enumerate}
\end{corollary}

Now, we obtain some classifications for the rotational hypersurface
(\ref{surftime}) which has $L_{1}$-harmonic Gauss map, first kind $L_{1}%
$-pointwise 1-type Gauss map, second kind $L_{1}$-pointwise 1-type Gauss map
and generalized $L_{1}$ 1-type Gauss map.

If $\Gamma^{t}(s,t,w)$ has a generalized $L_{1}$ (Cheng-Yau) 1-type Gauss map,
i.e., $L_{1}N^{t}=\mathfrak{m}N^{t}+\mathfrak{n}C,$ where $C=\left(
C_{1},C_{2},C_{3},C_{4}\right)  $ is constant vector, then from
(\ref{normaltime}) and (\ref{L1time}) we get%

\begin{equation}
\left.
\begin{array}
[c]{l}%
\frac{2R_{1}^{t}}{s^{3}\left(  f^{\prime2}-1\right)  ^{4}}=\mathfrak{m}\left(
\frac{1}{\sqrt{f^{\prime2}-1}}\right)  +\mathfrak{n}C_{1},\\
\\
\frac{2\cos t\sin wQ_{1}^{t}}{s^{3}\left(  f^{\prime2}-1\right)  ^{4}%
}=\mathfrak{m}\left(  \frac{-f^{\prime}\cos t\sin w}{\sqrt{f^{\prime2}-1}%
}\right)  +\mathfrak{n}C_{2},\\
\\
\frac{2\sin tQ_{1}^{t}}{s^{3}\left(  f^{\prime2}-1\right)  ^{4}}%
=\mathfrak{m}\left(  \frac{-f^{\prime}\sin t}{\sqrt{f^{\prime2}-1}}\right)
+\mathfrak{n}C_{3},\\
\\
\frac{-2\cos t\cos wQ_{1}^{t}}{s^{3}\left(  f^{\prime2}-1\right)  ^{4}%
}=\mathfrak{m}\left(  \frac{f^{\prime}\cos t\cos w}{\sqrt{f^{\prime2}-1}%
}\right)  +\mathfrak{n}C_{4}.
\end{array}
\right\}  \label{sinif1time}%
\end{equation}
It is obvious from the first three equations of (\ref{sinif1time}), we get
$C_{2}=C_{3}=C_{4}=0$. Here if $f^{\prime}=0,$ then we have $L_{1}N^{t}=0.$ Thus,

\begin{theorem}
The rotational hypersurface
\[
\Gamma^{t}(s,t,w)=\left(  c_{9},-s\cos t\sin w,-s\sin t,s\cos t\cos w\right)
\]
has $L_{1}$-harmonic Gauss map in $\mathbb E_{1}^{4}.$
\end{theorem}

From now on, we assume that $f^{\prime}\neq0.$ In (\ref{sinif1time}), we can
take $C_{1}=0$ or $C_{1}\neq0.$

Firstly let us suppose that $C_{1}=0.$ In this case; from the first equation
and the last three equations of (\ref{sinif1time}), we have%
\begin{equation}
\mathfrak{m}=\frac{2f^{\prime}\left(  \left(  f^{\prime}\left(  f^{\prime
2}-1\right)  -s\left(  3f^{\prime2}+2\right)  f^{\prime\prime}\right)
f^{\prime\prime}+s\left(  f^{\prime2}-1\right)  f^{\prime}f^{\prime
\prime\prime}\right)  }{s^{2}\left(  f^{\prime2}-1\right)  ^{7/2}}
\label{sinif2time}%
\end{equation}
and%
\begin{equation}
\mathfrak{m}=\frac{-2f^{\prime2}\left(  f^{\prime2}-1\right)  ^{3}+2s\left(
f^{\prime2}-1\right)  f^{\prime3}f^{\prime\prime}-2s^{2}\left(  4f^{\prime
2}+1\right)  f^{\prime\prime2}+2s^{2}\left(  f^{\prime2}-1\right)  f^{\prime
}f^{\prime\prime\prime}}{s^{3}f^{\prime}\left(  f^{\prime2}-1\right)  ^{7/2}},
\label{sinif3time}%
\end{equation}
respectively. From (\ref{sinif2time}) and (\ref{sinif3time}), we get%
\begin{equation}
\left(  f^{\prime}\left(  \left(  f^{\prime2}-1\right)  ^{2}-3s^{2}%
f^{\prime\prime2}\right)  +s^{2}\left(  f^{\prime2}-1\right)  f^{\prime
\prime\prime}\right)  f^{\prime}-s^{2}f^{\prime\prime2}=0. \label{sinif4time}%
\end{equation}
By solving (\ref{sinif4time}), we obtain that%
\[
f=\pm\int_{1}^{s}\frac{\sqrt{2\left(  6c_{10}c_{11}x^{3}+c_{10}\right)  }%
}{\sqrt{12c_{10}c_{11}x^{3}+3x+2c_{10}}}dx+c_{12}.
\]
Hence, we obtain the following result:
\begin{theorem}
The rotational hypersurface
\[
\Gamma^{t}(s,t,w)=\left(  \pm\int_{1}^{s}\frac{\sqrt{2\left(  6c_{10}%
c_{11}x^{3}+c_{10}\right)  }}{\sqrt{12c_{10}c_{11}x^{3}+3x+2c_{10}}}%
dx+c_{12},-s\cos t\sin w,-s\sin t,s\cos t\cos w\right)
\]
has first kind $L_{1}$-pointwise 1-type Gauss map, i.e., $L_{1}N^{t}%
=\mathfrak{m}N^{t},$ in $\mathbb E_{1}^{4},$ where%
\begin{align*}
&  \mathfrak{m}=\frac{c_{10}{}^{2}\sqrt{2}\left(  1+12c_{11}s^{3}\left(
1+12c_{11}s^{3}\right)  \right)  }{(-s)^{9/2}\sqrt{3c_{10}\left(
1+6c_{11}s^{3}\right)  }},\\
&  \text{and}\\
&  L_{1}N^{t}=\frac{2c_{10}{}^{2}\left(  1+12s^{3}c_{11}\left(  1+12c_{11}%
s^{3}\right)  \right)  }{3s^{5}}\left( \frac{\sqrt{6s+4c_{10}+24c_{10}%
c_{11}s^{3}}}{2\sqrt{c_{10}\left(  1+6c_{11}s^{3}\right)  }},\cos t\sin w,\sin
t,-\cos t\cos w\right)  .
\end{align*}
\end{theorem}

Now, we suppose that $C_{1}\neq0.$ In this case, from the last three
equations and the first equation of (\ref{sinif1time}), we have%
\begin{equation}
\left.
\begin{array}
[c]{l}%
\mathfrak{m}=-\frac{2\left(  3f^{\prime4}\left(  1-f^{\prime2}\right)
+f^{\prime8}-sf^{\prime5}f^{\prime\prime}+s^{2}f^{\prime\prime2}+\left(
4s^{2}f^{\prime\prime2}-1\right)  f^{\prime2}+s^{2}f^{\prime}f^{\prime
\prime\prime}+sf^{\prime3}\left(  f^{\prime\prime}-sf^{\prime\prime\prime
}\right)  \right)  }{s^{3}f^{\prime}\left(  f^{\prime2}-1\right)  ^{7/2}}\\
\text{and}\\
\mathfrak{n}=\frac{2\left(  f^{\prime}\left(  \left(  f^{\prime2}-1\right)
^{2}-3s^{2}f^{\prime\prime2}\right)  +s^{2}\left(  f^{\prime2}-1\right)
f^{\prime\prime\prime}\right)  f^{\prime}-2s^{2}f^{\prime\prime2}}{C_{1}%
s^{3}f^{\prime}\left(  f^{\prime2}-1\right)  ^{3}},
\end{array}
\right\}  \label{sinif5time}%
\end{equation}
respectively.

\begin{theorem}
The rotational hypersurface (\ref{surftime}) has generalized $L_{1}$
1-type Gauss map, i.e., $L_{1}N^{t}=\mathfrak{m}N^{t}+\mathfrak{n}C,$ in
$\mathbb E_{1}^{4},$ where $\mathfrak{m}$ and $\mathfrak{n}$ are non-zero smooth
functions given by (\ref{sinif5time}) and $C=(C_{1},0,0,0)$ is a non-zero
vector with non-zero constant $C_{1}$.
\end{theorem}

If we consider the case $\mathfrak{m}=\mathfrak{n}$ in (\ref{sinif1time}) (or in
(\ref{sinif5time})), then we get
\begin{align}\label{sinif7time}
\begin{split}
0=&  \left(  1+C_{1}\sqrt{f^{\prime2}-1}\right)  \left(
-2f^{\prime2}\left(  f^{\prime2}-1\right)  ^{3}+2s\left(  f^{\prime
2}-1\right)  f^{\prime3}f^{\prime\prime}\right.\\&\left. 
-2s^{2}\left(  4f^{\prime2}+1\right)  f^{\prime\prime2}+2s^{2}\left(
f^{\prime2}-1\right)  f^{\prime}f^{\prime\prime\prime}%
/\right)\\
& -2s\left(  \left(  f^{\prime}\left(  f^{\prime2}-1\right)  -s\left(
3f^{\prime2}+2\right)  f^{\prime\prime}\right)  f^{\prime\prime}+s\left(
f^{\prime2}-1\right)  f^{\prime}f^{\prime\prime\prime}\right)  f^{\prime2}.
\end{split}
\end{align}
Therefore,

\begin{theorem}
The rotational hypersurface (\ref{surftime}) has second kind $L_{1}$-pointwise
1-type Gauss map, i.e., $L_{1}N^{t}=\mathfrak{m}\left(  N^{t}+C\right)  ,$ in
$\mathbb E_{1}^{4},$ if the differential equation (\ref{sinif7time}) holds.
\end{theorem}

\subsection{Rotational Hypersurfaces about Lightlike Axis}

\

The functions $a_{k}$ of rotational hypersurface (\ref{surflight}) in
$\mathbb E_{1}^{4}$ are obtained from (\ref{ak}) and (\ref{aslilight}) by%
\begin{equation}
\left.
\begin{array}
[c]{l}%
a_{1}^{l}=\frac{2\left(  f^{\prime}-1\right)  ^{2}\left(  f^{\prime}+1\right)
-\left(  s-f\right)  f^{\prime\prime}}{\left(  1-f^{\prime2}\right)
^{3/2}\left(  s-f\right)  },\\
\\
a_{2}^{l}=\frac{-\left(  f^{\prime}-1\right)  ^{2}\left(  f^{\prime}+1\right)
+2(s-f)f^{\prime\prime}}{(s-f)^{2}\left(  f^{\prime}-1\right)  \left(
f^{\prime}+1\right)  ^{2}},\\
\\
a_{3}^{l}=-\frac{f^{\prime\prime}}{(s-f)^{2}\sqrt{1-f^{\prime}}\left(
f^{\prime}+1\right)  ^{5/2}}.
\end{array}
\right\}  \label{aklight}%
\end{equation}
Also, from (\ref{gradf}), (\ref{glight}) and (\ref{aklight}), we have%
\begin{align}\label{a2nablalight}
\begin{split}
\nabla a_{2}^{l}=\frac{-P_{1}^{l}\left(  f^{\prime}+1\right)  ^{-4}}%
{(s-f)^{3}\left(  f^{\prime}-1\right)  ^{3}}&\left(
\left(  t^{2}+w^{2}\right)  f^{\prime}-2-t^{2}-w^{2},\left(  -2+t^{2}%
+w^{2}\right)  f^{\prime}-t^{2}-w^{2},\right.\\&\left.
2t\left(  f^{\prime}-1\right)  ,2w\left(  f^{\prime}-1\right)
\right)  , %
\end{split}
\end{align}
where
\begin{align}\nonumber
\begin{split}
P_{1}^{l}=&(s-f)^{2}f^{\prime\prime2}+f^{\prime}\left(  2\left(  1-f^{\prime
2}\left(  2-f^{\prime2}\right)  \right)  +f^{\prime}\left(  1+f^{\prime
2}\right)  -3(s-f)^{2}f^{\prime\prime2}-f^{\prime5}\right)\\&  +(s-f)^{2}\left(
f^{\prime2}-1\right)  f^{\prime\prime\prime}-1.
\end{split}
\end{align}

So, from (\ref{Lk}), (\ref{Hk}), (\ref{normallight}), (\ref{epsilonlight}),
(\ref{aklight}) and (\ref{a2nablalight}), we have%
\begin{equation}
L_{1}N^{l}=\frac{-\left(  f^{\prime}+1\right)  ^{-4}}{(s-f)^{3}\left(
f^{\prime}-1\right)  ^{3}}\left(
\begin{array}
[c]{l}%
\left(  t^{2}+w^{2}\right)  \left(  f^{\prime}P_{1}^{l}-Q_{1}^{l}\right)
+\left(  2+t^{2}+w^{2}\right)  \left(  f^{\prime}Q_{1}^{l}-P_{1}^{l}\right)
,\\
\left(  t^{2}+w^{2}\right)  \left(  f^{\prime}Q_{1}^{l}-P_{1}^{l}\right)
+\left(  2-t^{2}-w^{2}\right)  \left(  Q_{1}^{l}-f^{\prime}P_{1}^{l}\right)
,\\
-2t(s-f)\left(  f^{\prime}-1\right)  R_{1}^{l},\\
-2w(s-f)\left(  f^{\prime}-1\right)  R_{1}^{l}%
\end{array}
\right)  , \label{L1light}%
\end{equation}
where
\begin{align*}
Q_{1}^{l}  &  =\left(  f^{\prime}-1\right)  ^{4}\left(  f^{\prime}+1\right)
^{2}-(s-f)\left(  \left(  f^{\prime}-1\right)  ^{2}\left(  f^{\prime
}+1\right)  f^{\prime\prime}-(s-f)f^{\prime\prime2}\right)  ,\\
R_{1}^{l}  &  =\left(  f^{\prime}-1\right)  ^{2}\left(  f^{\prime}+1\right)
f^{\prime\prime}+(s-f)\left(  \left(  3f^{\prime}-2\right)  f^{\prime\prime
2}-\left(  f^{\prime2}-1\right)  f^{\prime\prime\prime}\right)  .
\end{align*}

If the rotational hypersurface (\ref{surflight}) in $\mathbb E_{1}^{4}$ is flat (resp.
minimal), then $a_{3}^{l}=0$ (resp. $a_{1}^{l}=0$)$.$ So, from (\ref{L1light}%
), we get

\begin{proposition}
If the rotational hypersurface (\ref{surflight}) in $\mathbb E_{1}^{4}$ is flat (resp.
minimal), then we obtain%
\begin{align}
&  L_{1}N^{l}=\frac{-\left(  f^{\prime}+1\right)  ^{-4}}{2(s-f)^{3}\left(
f^{\prime}-1\right)  ^{3}}\left(
\begin{array}
[c]{l}%
\left(  t^{2}+w^{2}+2\right)  \left(  f^{\prime}A_{1}^{l}-2P_{1}^{l}\right)
-\left(  t^{2}+w^{2}\right)  \left(  A_{1}^{l}-2f^{\prime}P_{1}^{l}\right)
,\\
\left(  t^{2}+w^{2}-2\right)  \left(  2f^{\prime}P_{1}^{l}-A_{1}^{l}\right)
+\left(  t^{2}+w^{2}\right)  \left(  f^{\prime}A_{1}^{l}-2P_{1}^{l}\right)
,\\
2t(s-f)\left(  1-f^{\prime}\right)  A_{2}^{l},\\
2w(s-f)\left(  1-f^{\prime}\right)  A_{2}^{l}%
\end{array}
\right) \label{L1lightflat}\\
&  \left(  \text{resp. }L_{1}N^{l}=\frac{-\left(  f^{\prime}+1\right)  ^{-4}%
}{2(s-f)^{3}\left(  f^{\prime}-1\right)  ^{3}}\left(
\begin{array}
[c]{l}%
\left(  t^{2}+w^{2}+2\right)  \left(  -f^{\prime}A_{4}^{l}-2A_{3}^{l}\right)
+\left(  t^{2}+w^{2}\right)  \left(  A_{4}^{l}+2f^{\prime}A_{3}^{l}\right)
,\\
\left(  t^{2}+w^{2}-2\right)  \left(  2f^{\prime}A_{3}^{l}+A_{4}^{l}\right)
-\left(  t^{2}+w^{2}\right)  \left(  2A_{3}^{l}+f^{\prime}A_{4}^{l}\right)
,\\
2t\left(  1-f^{\prime}\right)  (A_{4}^{l}-2A_{3}^{l}),\\
2w\left(  1-f^{\prime}\right)  (A_{4}^{l}-2A_{3}^{l})
\end{array}
\right)  \right)  , \label{L1lightminimal}%
\end{align}
where
\begin{align*}
&  A_{1}^{l}=\left(  (s-f)f^{\prime\prime}-2\left(  f^{\prime}-1\right)
^{2}\left(  f^{\prime}+1\right)  \right)  \left(  2(s-f)f^{\prime\prime
}-\left(  f^{\prime}-1\right)  ^{2}\left(  f^{\prime}+1\right)  \right)  ,\\
&  A_{2}^{l}=5\left(  f^{\prime}+1\right)  \left(  f^{\prime}-1\right)
^{2}f^{\prime\prime}-2(s-f)\left(  f^{\prime2}-1\right)  f^{\prime\prime
\prime}+2(s-f)\left(  3f^{\prime}-2\right)  f^{\prime\prime2},\\
&  \left(
\begin{array}
[c]{l}%
\text{resp. }A_{3}^{l}=(s-f)^{2}\left(  f^{\prime\prime2}\left(  1-3f^{\prime
}\right)  +\left(  f^{\prime2}-1\right)  f^{\prime\prime\prime}\right)
+f^{\prime}\left(  f^{\prime3}\left(  1-f^{\prime2}\right)  -2f^{\prime
2}\left(  2-f^{\prime2}\right)  +f^{\prime}+2\right)  -1,\\
A_{4}^{l}=-3(s-f)\left(  f^{\prime}-1\right)  ^{2}\left(  f^{\prime}+1\right)
f^{\prime\prime}%
\end{array}
\right)  .
\end{align*}

\end{proposition}

\begin{remark}
\label{remarklight}We know that \cite{hacettepe}, if the rotational
hypersurface (\ref{surflight}) in $\mathbb E_{1}^{4}$ is flat, then we have
$f(s)=c_{13}s+c_{14}$.
\end{remark}

Thus, from (\ref{L1light}) or (\ref{L1lightflat}) (resp. (\ref{L1lightminimal}%
)), we have

\begin{corollary}
If the rotational hypersurface (\ref{surflight}) in $\mathbb E_{1}^{4}$ is flat, then we
have%
\[
L_{1}N^{l}=\frac1{\left(  c_{13}+1\right)\left(  c_{14}+\left(  c_{13}-1\right)  s\right)  ^{3}}
\left(  2\left(  c_{13}-1\right)  ,2\left(c_{13}-1\right)  ,0,0\right)  .
\]

\end{corollary}

\section{On Hypersurfaces with Generalized $L_2$ 1-Type Gauss Map}

In this section, we obtain the $L_{2}$ operator of the Gauss maps of the
rotational hypersurfaces (\ref{surfspace}), (\ref{surftime}), (\ref{surflight}%
) and give some classifications for these hypersurfaces which have $L_{2}%
$-harmonic Gauss map, first kind $L_{2}$-pointwise 1-type Gauss map,
second kind $L_{2}$-pointwise 1-type Gauss map and generalized $L_{2}%
$ 1-type Gauss map, seperately.

\subsection{ Rotational Hypersurfaces about Spacelike Axis}

\

The gradient of function $a_{3}$ of the rotational hypersurface
(\ref{surfspace}) in $\mathbb E_{1}^{4}$ are obtained from (\ref{gradf}),
(\ref{gspace}) and (\ref{akspace}) by%

\begin{equation}
\nabla a_{3}^{s}=\frac{P_{2}^{s}}{s^{3}\left(  1-f^{\prime2}\right)  ^{9/2}%
}\left(  \cosh t\cosh w,\sinh t,\cosh t\sinh w,f^{\prime}\right)  ,
\label{a3nablaspace}%
\end{equation}
where
\[
P_{2}^{s}=f^{\prime}\left(  \left(  2f^{\prime}\left(  f^{\prime2}-1\right)
+s\left(  3f^{\prime2}+2\right)  f^{\prime\prime}\right)  f^{\prime\prime
}-s\left(  f^{\prime2}-1\right)  f^{\prime}f^{\prime\prime\prime}\right)  .
\]

So, from (\ref{Lk}), (\ref{Hk}), (\ref{normalspace}), (\ref{epsilonspace}),
(\ref{akspace}) and (\ref{a3nablaspace}), we have%
\begin{equation}
L_{2}N^{s}=\frac{1}{s^{3}\left(  1-f^{\prime2}\right)  ^{9/2}}\left(  \cosh
t\cosh wQ_{2}^{s},\sinh tQ_{2}^{s},\cosh t\sinh wQ_{2}^{s},R_{2}^{s}\right)
, \label{L2space}%
\end{equation}
where
\begin{align*}
Q_{2}^{s}  &  =f^{\prime}\left(  2\left(  -f^{\prime}\left(  f^{\prime
2}-1\right)  ^{2}+s\left(  2f^{\prime2}+1\right)  f^{\prime\prime}\right)
f^{\prime\prime}-s\left(  f^{\prime2}-1\right)  f^{\prime}f^{\prime
\prime\prime}\right)  ,\\
R_{2}^{s}  &  =sf^{\prime2}\left(  3\left(  f^{\prime2}+1\right)
f^{\prime\prime2}-\left(  f^{\prime2}-1\right)  f^{\prime}f^{\prime
\prime\prime}\right)
\end{align*}

\begin{proposition}
If the rotational hypersurface (\ref{surfspace}) in $\mathbb E_{1}^{4}$ is flat (resp.
minimal), then we obtain%
\begin{align}
&  L_{2}N^{s}=0\label{L2spaceflat}\\
&  \left(  \text{resp. }L_{2}N^{s}=\frac{A_{5}^{s}}{s^{3}\left(  1-f^{\prime
2}\right)  ^{9/2}}\left(
\begin{array}
[c]{c}%
\cosh t\cosh w,\sinh t,\cosh t\sinh w,f^{\prime}%
\end{array}
\right)  \right)  , \label{L2spaceminimal}%
\end{align}
where
\[
A_{5}^{s}=f^{\prime}\left(  \left(  2f^{\prime}\left(  f^{\prime2}-1\right)
+s\left(  3f^{\prime2}+2\right)  f^{\prime\prime}\right)  f^{\prime\prime
}-s\left(  f^{\prime2}-1\right)  f^{\prime}f^{\prime\prime\prime}\right)  .
\]

\end{proposition}

Thus, using Remark \ref{remarkspace}, from (\ref{L2space}) or
(\ref{L2spaceminimal}), we have

\begin{corollary}
If the rotational hypersurface (\ref{surfspace}) in $\mathbb E_{1}^{4}$ is minimal, then we
have%
\begin{align}\nonumber
\begin{split}
L_{2}N^{s}=\frac{18c_{2}{}^{3}}{s^{14}\sqrt{s^{4}+c_{2}{}^{2}}}&\left(
\left(  s^{4}+c_{2}{}^{2}\right)  {}^{3/2}\cosh t\cosh w,\left(  s^{4}+c_{2}%
{}^{2}\right)  {}^{3/2}\sinh t,\right.\\&\left.
\left(  s^{4}+c_{2}{}^{2}\right)  {}^{3/2}\cosh t\sinh w,c_{2}\left(
s^{4}+c_{2}{}^{2}\right)
\right).
\end{split}
\end{align}

\end{corollary}

Now, let us give some classifications for the rotational hypersurface
(\ref{surfspace}) which has $L_{2}$-harmonic Gauss map, first kind $L_{2}%
$-pointwise 1-type Gauss map, second kind $L_{2}$-pointwise 1-type Gauss map
and generalized $L_{2}$ 1-type Gauss map.

If $\Gamma^{s}(s,t,w)$ has a generalized $L_{2}$ 1-type Gauss map, i.e.,
$L_{2}N^{s}=\mathfrak{m}N^{s}+\mathfrak{n}C,$ where $C=\left(  C_{1}%
,C_{2},C_{3},C_{4}\right)  $ is constant vector, then from (\ref{normalspace})
and (\ref{L2space}) we get%
\begin{equation}
\left.
\begin{array}
[c]{l}%
\frac{\cosh t\cosh wQ_{2}^{s}}{s^{3}\left(  1-f^{\prime2}\right)  ^{9/2}%
}=\mathfrak{m}\left(  \frac{-f^{\prime}\cosh t\cosh w}{\sqrt{1-f^{\prime2}}%
}\right)  +\mathfrak{n}C_{1},\\
\\
\frac{\sinh tQ_{2}^{s}}{s^{3}\left(  1-f^{\prime2}\right)  ^{9/2}%
}=\mathfrak{m}\left(  \frac{-f^{\prime}\sinh t}{\sqrt{1-f^{\prime2}}}\right)
+\mathfrak{n}C_{2},\\
\\
\frac{\cosh t\sinh wQ_{2}^{s}}{s^{3}\left(  1-f^{\prime2}\right)  ^{9/2}%
}=\mathfrak{m}\left(  \frac{-f^{\prime}\cosh t\sinh w}{\sqrt{1-f^{\prime2}}%
}\right)  +\mathfrak{n}C_{3},\\
\\
\frac{R_{2}^{s}}{s^{3}\left(  1-f^{\prime2}\right)  ^{9/2}}=\mathfrak{m}%
\left(  \frac{-1}{\sqrt{1-f^{\prime2}}}\right)  +\mathfrak{n}C_{4}.
\end{array}
\right\}  \label{sinif1space2}%
\end{equation}
It is obvious from the first three equations of (\ref{sinif1space2}), we get
$C_{1}=C_{2}=C_{3}=0$. Here if $f^{\prime\prime}=0,$ then we have $L_{2}%
N^{s}=0.$ Thus,

\begin{theorem}
The rotational hypersurface
\[
\Gamma^{s}(s,t,w)=\left(  s\cosh t\cosh w,s\sinh t,s\cosh t\sinh
w,d_{1}s+d_{2}\right)  ,\text{ }\left(  d_{1}\in \mathbb R-\left\{  \pm1\right\}
,\text{ }d_{2}\in\mathbb R\right)
\]
has $L_{2}$-harmonic Gauss map in $\mathbb E_{1}^{4}.$
\end{theorem}

From now on, we assume that $f^{\prime\prime}\neq0.$ In (\ref{sinif1space2}),
we can take $C_{4}=0$ or $C_{4}\neq0.$

Firstly let us suppose that $C_{4}=0.$ In this case; from the first three
equations and the last equation of (\ref{sinif1space2}), we have%
\begin{equation}
\mathfrak{m}=\frac{s\left(  f^{\prime2}-1\right)  f^{\prime}f^{\prime
\prime\prime}+2\left(  f^{\prime}\left(  f^{\prime2}-1\right)  ^{2}-s\left(
2f^{\prime2}+1\right)  f^{\prime\prime}\right)  f^{\prime\prime}}{s^{3}\left(
f^{\prime2}-1\right)  ^{4}} \label{sinif2space2}%
\end{equation}
and%
\begin{equation}
\mathfrak{m}=\frac{f^{\prime2}\left(  \left(  f^{\prime2}-1\right)  f^{\prime
}f^{\prime\prime\prime}-3\left(  f^{\prime2}+1\right)  f^{\prime\prime
2}\right)  }{s^{2}\left(  f^{\prime2}-1\right)  ^{4}}, \label{sinif3space2}%
\end{equation}
respectively. From (\ref{sinif2space2}) and (\ref{sinif3space2}), we get%
\begin{equation}
s\left(  f^{\prime2}-1\right)  f^{\prime}f^{\prime\prime\prime}-\left(
2f^{\prime}\left(  f^{\prime2}-1\right)  +s\left(  3f^{\prime2}+2\right)
f^{\prime\prime}\right)  f^{\prime\prime}=0. \label{sinif4space2}%
\end{equation}
Hence, from (\ref{sinif4space2}), we have

\begin{theorem}
The rotational hypersurface (\ref{surfspace}) has first kind $L_{2}$-pointwise
1-type Gauss map, i.e., $L_{2}N^{s}=\mathfrak{m}N^{s},$ in $\mathbb E_{1}^{4},$ if the
differential equation%
\[
\frac{f^{\prime\prime\prime}}{f^{\prime\prime}}=\frac{s\left(  3f^{\prime
2}+2\right)  f^{\prime\prime}+2f^{\prime}\left(  f^{\prime2}-1\right)
}{sf^{\prime}\left(  1-f^{\prime2}\right)  }%
\]
holds.
\end{theorem}

Now let us suppose that $C_{4}\neq0.$ In this case; from the first three
equations and the last equation of (\ref{sinif1space2}), we have%
\begin{equation}
\left.
\begin{array}
[c]{l}%
\mathfrak{m}=\frac{2f^{\prime\prime}\left(  f^{\prime}\left(  f^{\prime
2}-1\right)  ^{2}-s\left(  2f^{\prime2}+1\right)  f^{\prime\prime}\right)
+s\left(  f^{\prime2}-1\right)  f^{\prime}f^{\prime\prime\prime}}{s^{3}\left(
f^{\prime2}-1\right)  ^{4}}\\
\text{and}\\
\mathfrak{n}=-\frac{f^{\prime\prime}\left(  2f^{\prime}\left(  f^{\prime
2}-1\right)  +s\left(  3f^{\prime2}+2\right)  f^{\prime\prime}\right)
-s\left(  f^{\prime2}-1\right)  f^{\prime}f^{\prime\prime\prime}}{C_{4}%
s^{3}\left(  1-f^{\prime2}\right)  ^{7/2}},
\end{array}
\right\}  \label{sinif5space2}%
\end{equation}
respectively.

\begin{theorem}
The rotational hypersurface (\ref{surfspace}) has generalized $L_{2}%
$ 1-type Gauss map, i.e., $L_{2}N^{s}=\mathfrak{m}N^{s}%
+\mathfrak{n}C,$ in $\mathbb E_{1}^{4},$ where $\mathfrak{m}$ and $\mathfrak{n}$ are non-zero smooth functions
given by (\ref{sinif5space2}) and $C=(0,0,0,C_{4})$ is a non-zero vector with
non-zero constant $C_{4}$.
\end{theorem}

If we take $\mathfrak{m}=\mathfrak{n}$ in (\ref{sinif1space2}) (or in
(\ref{sinif5space2})), then we obtain%
\begin{align}
&  \left(  C_{4}\sqrt{1-f^{\prime2}}-1\right)  \left(  s\left(  f^{\prime
2}-1\right)  f^{\prime}f^{\prime\prime\prime}+2f^{\prime\prime}\left(
f^{\prime}\left(  f^{\prime2}-1\right)  ^{2}-s\left(  2f^{\prime2}+1\right)
f^{\prime\prime}\right)  \right) \nonumber\\
&  +sf^{\prime2}\left(  \left(  f^{\prime2}-1\right)  f^{\prime}%
f^{\prime\prime\prime}-3\left(  f^{\prime2}+1\right)  f^{\prime\prime
2}\right)  =0. \label{sinif7space2}%
\end{align}
Therefore,

\begin{theorem}
The rotational hypersurface (\ref{surfspace}) has second kind $L_{2}%
$-pointwise 1-type Gauss map, i.e., $L_{2}N^{s}=\mathfrak{m}\left(
N^{s}+C\right)  ,$ in $\mathbb E_{1}^{4},$ if the differential equation
(\ref{sinif7space2}) holds.
\end{theorem}

\subsection{Rotational Hypersurfaces about Timelike Axis}

\

The gradient of function $a_{3}$ of the rotational hypersurface
(\ref{surftime}) in $\mathbb E_{1}^{4}$ are obtained from (\ref{gradf}), (\ref{gtime})
and (\ref{aktime}) by%

\begin{equation}
\nabla a_{3}^{t}=\frac{P_{2}^{t}}{s^{3}\left(  f^{\prime2}-1\right)  ^{9/2}%
}\left(  -f^{\prime},\cos t\sin w,\sin t,-\cos t\cos w\right)  ,
\label{a3nablatime}%
\end{equation}
where
\[
P_{2}^{t}=f^{\prime}\left(  \left(  2f^{\prime}\left(  f^{\prime2}-1\right)
+s\left(  3f^{\prime2}+2\right)  f^{\prime\prime}\right)  f^{\prime\prime
}-s\left(  f^{\prime2}-1\right)  f^{\prime}f^{\prime\prime\prime}\right)  .
\]

So, from (\ref{Lk}), (\ref{Hk}), (\ref{normaltime}), (\ref{epsilontime}),
(\ref{aktime}) and (\ref{a3nablatime}), we have%
\begin{equation}
L_{2}N^{t}=\frac{1}{s^{3}\left(  f^{\prime2}-1\right)  ^{9/2}}\left(
R_{2}^{t},\cos t\sin wQ_{2}^{t},\sin tQ_{2}^{t},-\cos t\cos wQ_{2}%
^{t}\right)  , \label{L2time}%
\end{equation}
where
\begin{align*}
Q_{2}^{t}  &  =f^{\prime}\left(  2\left(  s\left(  2f^{\prime2}+1\right)
f^{\prime\prime}-f^{\prime}\left(  f^{\prime2}-1\right)  ^{2}\right)
f^{\prime\prime}-s\left(  f^{\prime2}-1\right)  f^{\prime}f^{\prime
\prime\prime}\right)  ,\\
R_{2}^{t}  &  =sf^{\prime2}\left(  \left(  f^{\prime2}-1\right)  f^{\prime
}f^{\prime\prime\prime}-3\left(  f^{\prime2}+1\right)  f^{\prime\prime
2}\right)
\end{align*}

\begin{proposition}
If the rotational hypersurface (\ref{surftime}) in $\mathbb E_{1}^{4}$ is flat (resp.
minimal), then we obtain%
\begin{align}
&  L_{2}N^{t}=0\label{L2timeflat}\\
&  \left(  \text{resp. }L_{2}N^{t}=\frac{A_{5}^{t}}{s^{3}\left(  f^{\prime
2}-1\right)  ^{9/2}}\left(  -f^{\prime},\cos t\sin w,\sin t,-\cos t\cos
w\right)  \right)  , \label{L2timeminimal}%
\end{align}
where
\[
A_{5}^{t}=f^{\prime}\left(  f^{\prime\prime}\left(  2f^{\prime}\left(
f^{\prime2}-1\right)  +s\left(  3f^{\prime2}+2\right)  f^{\prime\prime
}\right)  -s\left(  f^{\prime2}-1\right)  f^{\prime}f^{\prime\prime\prime
}\right)  .
\]

\end{proposition}

Thus, using Remark \ref{remarktime}, from (\ref{L2time}) or
(\ref{L2timeminimal}), we have

\begin{corollary}
If the rotational hypersurface (\ref{surftime}) in $\mathbb E_{1}^{4}$ is minimal, we
have%
\begin{align}\nonumber
\begin{split}
L_{2}N^{t}=\frac{18c_{8}{}^{3}}{s^{14}\sqrt{c_{8}{}^{2}-s^{4}}}&\left(
c_{8}{}\left(  s^{4}-c_{8}{}^{2}\right)  ,\left(  c_{8}{}^{2}-s^{4}\right)
{}^{3/2}\cos t\sin w,\right.\\&\left.
\left(  c_{8}{}^{2}-s^{4}\right)  {}^{3/2}\sin t,-\left(  c_{8}{}^{2}%
-s^{4}\right)  {}^{3/2}\cos t\cos w
\right)  .
\end{split}
\end{align}

\end{corollary}

Now, let us give some classifications for the rotational hypersurface
(\ref{surftime}) which has $L_{2}$-harmonic Gauss map, first kind $L_{2}%
$-pointwise 1-type Gauss map, second kind $L_{2}$-pointwise 1-type Gauss map
and generalized $L_{2}$ 1-type Gauss map.

If $\Gamma^{t}(s,t,w)$ has a generalized $L_{2}$ 1-type Gauss map, i.e.,
$L_{2}N^{t}=\mathfrak{m}N^{t}+\mathfrak{n}C,$ where $C=\left(  C_{1}%
,C_{2},C_{3},C_{4}\right)  $ is constant vector, then from (\ref{normaltime})
and (\ref{L2time}) we get%

\begin{equation}
\left.
\begin{array}
[c]{l}%
\frac{R_{2}^{t}}{s^{3}\left(  f^{\prime2}-1\right)  ^{9/2}}=\mathfrak{m}%
\left(  \frac{1}{\sqrt{f^{\prime2}-1}}\right)  +\mathfrak{n}C_{1},\\
\\
\frac{\cos t\sin wQ_{2}^{t}}{s^{3}\left(  f^{\prime2}-1\right)  ^{9/2}%
}=\mathfrak{m}\left(  \frac{-f^{\prime}\cos t\sin w}{\sqrt{f^{\prime2}-1}%
}\right)  +\mathfrak{n}C_{2},\\
\\
\frac{\sin tQ_{2}^{t}}{s^{3}\left(  f^{\prime2}-1\right)  ^{9/2}}%
=\mathfrak{m}\left(  \frac{-f^{\prime}\sin t}{\sqrt{f^{\prime2}-1}}\right)
+\mathfrak{n}C_{3},\\
\\
\frac{-\cos t\cos wQ_{2}^{t}}{s^{3}\left(  f^{\prime2}-1\right)  ^{9/2}%
}=\mathfrak{m}\left(  \frac{f^{\prime}\cos t\cos w}{\sqrt{f^{\prime2}-1}%
}\right)  +\mathfrak{n}C_{4}.
\end{array}
\right\}  \label{sinif1time2}%
\end{equation}
It is obvious from the first three equations of (\ref{sinif1time2}), we get
$C_{2}=C_{3}=C_{4}=0$. Here if $f^{\prime\prime}=0,$ then we have $L_{2}%
N^{t}=0.$ Thus,

\begin{theorem}
The rotational hypersurface
\[
\Gamma^{t}(s,t,w)=\left(  d_{3}s+d_{4},-s\cos t\sin w,-s\sin t,s\cos t\cos
w\right)  ,\text{ }\left(  d_{3}\in \mathbb R-\{\pm1\},\text{ }d_{4}\in\mathbb  R\right)
\]
has $L_{2}$-harmonic Gauss map in $\mathbb E_{1}^{4}.$
\end{theorem}

From now on, we assume that $f^{\prime\prime}\neq0.$ In (\ref{sinif1time2}),
we can take $C_{1}=0$ or $C_{1}\neq0.$

Firstly let us suppose that $C_{1}=0.$ In this case; from the first equation
and the last three equations of (\ref{sinif1time2}), we have%
\begin{equation}
\mathfrak{m}=\frac{f^{\prime2}\left(  \left(  f^{\prime2}-1\right)  f^{\prime
}f^{\prime\prime\prime}-3\left(  f^{\prime2}+1\right)  f^{\prime\prime
2}\right)  }{s^{2}\left(  f^{\prime2}-1\right)  ^{4}} \label{sinif2time2}%
\end{equation}
and%
\begin{equation}
\mathfrak{m}=\frac{2f^{\prime\prime}\left(  f^{\prime}\left(  f^{\prime
2}-1\right)  ^{2}-s\left(  2f^{\prime2}+1\right)  f^{\prime\prime}\right)
+s\left(  f^{\prime2}-1\right)  f^{\prime}f^{\prime\prime\prime}}{s^{3}\left(
f^{\prime2}-1\right)  ^{4}}, \label{sinif3time2}%
\end{equation}
respectively. From (\ref{sinif2time2}) and (\ref{sinif3time2}), we get%
\begin{equation}
s\left(  f^{\prime2}-1\right)  f^{\prime}f^{\prime\prime\prime}-\left(
2f^{\prime}\left(  f^{\prime2}-1\right)  +s\left(  3f^{\prime2}+2\right)
f^{\prime\prime}\right)  f^{\prime\prime}=0. \label{sinif4time2}%
\end{equation}
Hence, from (\ref{sinif4time2}), we have

\begin{theorem}
The rotational hypersurface (\ref{surftime}) has first kind $L_{2}$-pointwise
1-type Gauss map, i.e., $L_{2}N^{t}=\mathfrak{m}N^{t},$ in $\mathbb E_{1}^{4},$ if the
differential equation%
\[
\frac{f^{\prime\prime\prime}}{f^{\prime\prime}}=\frac{s\left(  3f^{\prime
2}+2\right)  f^{\prime\prime}+2f^{\prime}\left(  f^{\prime2}-1\right)
}{sf^{\prime}\left(  1-f^{\prime2}\right)  }%
\]
holds.
\end{theorem}

Now let us suppose that $C_{1}\neq0.$ In this case; from the last three
equations and the first equation of (\ref{sinif1time2}), we have%
\begin{equation}
\left.
\begin{array}
[c]{l}%
\mathfrak{m}=\frac{2f^{\prime\prime}\left(  f^{\prime}\left(  f^{\prime
2}-1\right)  ^{2}-s\left(  2f^{\prime2}+1\right)  f^{\prime\prime}\right)
+s\left(  f^{\prime2}-1\right)  f^{\prime}f^{\prime\prime\prime}}{s^{3}\left(
f^{\prime2}-1\right)  ^{4}}\\
\text{and}\\
\mathfrak{n}=\frac{s\left(  f^{\prime2}-1\right)  f^{\prime}f^{\prime
\prime\prime}-\left(  2f^{\prime}\left(  f^{\prime2}-1\right)  +s\left(
3f^{\prime2}+2\right)  f^{\prime\prime}\right)  f^{\prime\prime}}{C_{1}%
s^{3}\left(  f^{\prime2}-1\right)  ^{7/2}},
\end{array}
\right\}  \label{sinif5time2}%
\end{equation}
respectively.

\begin{theorem}
The rotational hypersurface (\ref{surftime}) has generalized $L_{2}$
1-type Gauss map, i.e., $L_{2}N^{t}=\mathfrak{m}N^{t}+\mathfrak{n}C,$ in
$\mathbb E_{1}^{4},$ where $\mathfrak{m}$ and $\mathfrak{n}$ are non-zero smooth functions given by
(\ref{sinif5time2}) and $C=(C_{1},0,0,0)$ is a non-zero vector with non-zero
constant $C_{1}$.
\end{theorem}

If we take $\mathfrak{m}=\mathfrak{n}$ in (\ref{sinif1time2}) (or in
(\ref{sinif5time2})), then we obtain%
\begin{align}
&  \left(  1+C_{1}\sqrt{f^{\prime2}-1}\right)  \left(  2\left(  f^{\prime
}\left(  f^{\prime2}-1\right)  ^{2}-s\left(  2f^{\prime2}+1\right)
f^{\prime\prime}\right)  f^{\prime\prime}+s\left(  f^{\prime2}-1\right)
f^{\prime}f^{\prime\prime\prime}\right) \nonumber\\
&  -sf^{\prime2}\left(  \left(  f^{\prime2}-1\right)  f^{\prime}%
f^{\prime\prime\prime}-3\left(  f^{\prime2}+1\right)  f^{\prime\prime
2}\right)  =0. \label{sinif7time2}%
\end{align}
Therefore,

\begin{theorem}
The rotational hypersurface (\ref{surftime}) has second kind $L_{2}$-pointwise
1-type Gauss map, i.e., $L_{2}N^{t}=\mathfrak{m}\left(  N^{t}+C\right)  ,$ in
$\mathbb E_{1}^{4},$ if the differential equation (\ref{sinif7time2}) holds.
\end{theorem}

\subsection{Rotational Hypersurfaces about Lightlike Axis}

\

From (\ref{gradf}), (\ref{glight}) and (\ref{aklight}), we have%
\begin{align}\label{a3nablalight}%
\begin{split}
\nabla a_{3}^{l}=\frac{P_{2}^{l}\left(  f^{\prime}+1\right)  ^{-9/2}%
}{2(s-f)^{3}\left(  1-f^{\prime}\right)  ^{5/2}}&\left(
\left(  \left(  t^{2}+w^{2}\right)  f^{\prime}-t^{2}-w^{2}-2\right)  ,\left(
\left(  t^{2}+w^{2}-2\right)  f^{\prime}-t^{2}-w^{2}\right)  ,\right.\\&\left.
2\left(  f^{\prime}-1\right)  t,2\left(  f^{\prime}-1\right)  w
\right) 
\end{split}
\end{align}
where
\[
P_{2}^{l}=(s-f)\left(  f^{\prime2}-1\right)  f^{\prime\prime\prime}+2\left(
f^{\prime}+1\right)  \left(  f^{\prime}-1\right)  ^{2}f^{\prime\prime
}-(s-f)\left(  3f^{\prime}-2\right)  f^{\prime\prime2}.
\]

So, from (\ref{Lk}), (\ref{Hk}), (\ref{normallight}), (\ref{epsilonlight}),
(\ref{aklight}) and (\ref{a3nablalight}), we have%
\begin{equation}
L_{2}N^{l}=\frac{\left(  f^{\prime}+1\right)  ^{-9/2}}{2(s-f)^{2}\left(
1-f^{\prime}\right)  ^{5/2}}\left(\mathcal F_1(u,t,w),\mathcal F_2(u,t,w),2t\mathcal F_3(u,t,w),2w\mathcal F_3(u,t,w)\right)  , \label{L2light}%
\end{equation}
where we put
\begin{eqnarray*}
\mathcal F_1(u,t,w)&=& \left(  f^{\prime2}-1\right)  \left(  \left(  t^{2}+w^{2}\right)  f^{\prime
}-t^{2}-w^{2}-2\right)  f^{\prime\prime\prime}\\&&
-\left(3\left(  t^{2}+w^{2}\right)  f^{\prime2}-\left(  3t^{2}+3w^{2}+4\right)  \left(  2f^{\prime}-1\right)
\right)  f^{\prime\prime2}-\frac{4\left(  f^{\prime2}-1\right)  ^{2}%
f^{\prime\prime}}{(s-f)},\\
\mathcal F_2(u,t,w)&=& \left(  f^{\prime2}-1\right)  \left(  \left(  t^{2}+w^{2}-2\right)  f^{\prime
}-t^{2}-w^{2}\right)  f^{\prime\prime\prime}\\&&
-\left(3\left(  t^{2}+w^{2}-2\right)  f^{\prime2}-\left(  3t^{2}+3w^{2}-2\right)  \left(  2f^{\prime}-1\right)
\right)  f^{\prime\prime2}-\frac{4\left(  f^{\prime2}-1\right)  ^{2}%
f^{\prime\prime}}{(s-f)},\\
\mathcal F_3(u,t,w)&=& 
\left(  1-f^{\prime}\right)  ^{2}\left(  \left(  f^{\prime}+1\right)
f^{\prime\prime\prime}-3f^{\prime\prime2}\right).
\end{eqnarray*}

If the rotational hypersurface (\ref{surflight}) in $\mathbb E_{1}^{4}$ is flat (resp.
minimal), then $a_{3}^{l}=0$ (resp. $a_{1}^{l}=0$)$.$ So, from (\ref{L2light}%
), we get

\begin{proposition}
Let $M$ be the rotational hypersurface given by (\ref{surflight}) in $\mathbb E_{1}^{4}$. Then, we have the followings:
\begin{enumerate}
\item[(i)] If $M$ is flat, then its Gauss map satisfies $L_{2}N^{l}=0$.

\item[(ii)] If $M$ is minimal, then its Gauss map satisfies
\begin{align}\label{L2minimalflat}%
\begin{split}
L_{2}N^{l}=\frac{A_{5}^{l}\left(  f^{\prime}+1\right)
^{-9/2}}{2(s-f)^{3}\left(  1-f^{\prime}\right)  ^{5/2}}&\left(
\left(  t^{2}+w^{2}\right)  f^{\prime}-t^{2}-w^{2}-2,\left(  t^{2}%
+w^{2}-2\right)  f^{\prime}-t^{2}-w^{2},\right.\\&\left.
-2t\left(  1-f^{\prime}\right)  ,-2w\left(  1-f^{\prime}\right)
\right),
\end{split}
\end{align}
\end{enumerate}
where
\[
A_{5}^{l}=(s-f)\left(  \left(  f^{\prime2}-1\right)  f^{\prime\prime\prime
}-\left(  3f^{\prime}-2\right)  f^{\prime\prime2}\right)  +2\left(  f^{\prime
}+1\right)  \left(  f^{\prime}-1\right)  ^{2}f^{\prime\prime}.
\]
\end{proposition}


\begin{thebibliography}{99}                                                                                               %


\bibitem {Alias}L.J. Al\'{\i}as and N. G\"{u}rb\"{u}z, \textit{An extension of
Takahashi theorem for the linearized operators of the higher order mean
curvatures}, Geom. Dedicata \textbf{121} (2006), no. 1, 113--127, DOI 10.1007/s10711-006-9093-9.

\bibitem {Alias2}L. Al\'{\i}as, A. Ferr\'{a}ndez, P. Lucas, \textit{Surfaces
in the 3-dimensional Lorentz-Minkowski space satisfying }$\Delta x=Ax+B$,
Pacific J. Math. \textbf{156} (1992), no.2, 201-208.

\bibitem {hacettepe}M. Alt\i n and A. Kazan, \textit{Rotational Hypersurfaces
in Lorentz-Minkowski 4-Space}, Hacet. J. Math. Stat. \textbf{50} (2021), no.
5, 1409--1433.

\bibitem {Chen}B.Y. Chen and M. Petrovic, \textit{On spectral decomposition of
immersions of finite type,} Bull. Aust. Math. Soc. \textbf{44} (1991), no. 1, 117-129.

\bibitem{Chen-Piccinni} { Chen, B.-Y. and P. Piccinni,} \emph{Submanifolds with Finite Type Gauss Map}, Bull. Austral. Math. Soc., \textbf{35} (1987), 161-186.

\bibitem {Cheng}S.Y. Cheng and S.T. Yau, \textit{Hypersurfaces with constant
scalar curvature}, Math. Ann. \textbf{225} (1977), 195--204.

\bibitem {Dillen}F. Dillen, J. Pas and L. Verstraelen, \textit{On surfaces of
finite type in euclidean 3-space}, Kodai Math. J. \textbf{13} (1990), no. 1, 10-21.

\bibitem {Garay}O.J Garay, \textit{An extension of Takahashi's theorem}, Geom.
Dedicata \textbf{34} (1990), 105-112.

\bibitem {Guler}E. G\"{u}ler and N.C. Turgay, \textit{Cheng--Yau Operator and
Gauss Map of Rotational Hypersurfaces in 4-Space,} Mediterr. J. Math.
\textbf{16} (2019), no. 3: 66.

\bibitem {Hasanis}T. Hasanis and T. Vlachos, \textit{Hypersurfaces of
}$\mathbb E^{n+1}$\textit{ satisfying }$\Delta x=Ax+B$, J. Aust. Math. Soc.
\textbf{53} (1992), no. 3, 377-384.

\bibitem {Kelleci}A. Kelleci, Rotational surfaces with Cheng-Yau operator in
Galilean 3-spaces, Hacet. J. Math. Stat. \textbf{50} (2021), no. 2, 365--376,
DOI : 10.15672/hujms.612730.

\bibitem {Kim2}D-S. Kim, J.R. Kim and Y.H. Kim, \textit{Cheng-Yau Operator and
Gauss Map of Surfaces of Revolution}, Bull. Malays. Math. Sci. Soc.
\textbf{39} (2016), no.4, 1319--1327.

\bibitem {Kim}Y.H. Kim and N.C. Turgay, \textit{On Pointwise 1-Type Gauss Map
of Surfaces in }$\mathbb E_{1}^{3}$\textit{ Concerning Cheng-Yau Operator}, J. Korean
Math. Soc. \textbf{54} (2017), no. 2, 381--397.

\bibitem {Lucas}P. Lucas and H.F. Ram\'{\i}rez-Ospina, \textit{Hypersurfaces
in non-flat pseudo-Riemannian space forms satisfying a linear condition in the
linearized operator of a higher order mean curvature}, Taiwanese J. Math.
\textbf{17} (2013), no. 1, 15-45.

\bibitem {Lucas2}P. Lucas and H.F. Ram\'{\i}rez-Ospina, \textit{Hypersurfaces
in non-flat Lorentzian space forms satisfying }$L_{k}\psi=A\psi+b$, Taiwanese
J. Math. \textbf{16} (2012), no. 3, 1173-1203.

\bibitem {ChengYauEn1}P. Lucas and H.F. Ram\'{\i}rez-Ospina,
\textit{Hypersurfaces in the Lorentz-Minkowski space satisfying }$L_{k}%
\psi=A\psi+b$, Geom. Dedicata \textbf{153} (2011), no. 1, 151-175.

\bibitem {Qian}J. Qian, X. Fu, X. Tian and Y.H. Kim, \textit{Surfaces of
Revolution and Canal Surfaces with Generalized Cheng-Yau 1-Type Gauss Map},
Mathematics \textbf{8} (2020), no. 10, 1728.

\bibitem {Reilly}R.C. Reilly, \textit{Variational properties of functions of
the mean curvatures for hypersurfaces in space forms}, J. Diff. Geom.
\textbf{8} (1973), no.3, 465--477.

\bibitem {Rosenberg}H. Rosenberg, \textit{Hypersurfaces of constant curvature
in space forms}, Bull. Sci. Math. \textbf{117} (1993), no. 2, 211--239.

\bibitem {Takahashi}T. Takahashi, \textit{Minimal immersions of Riemannian
manifolds}, J. Math. Soc. Japan \textbf{18} (1966), no. 4, 380-385.

\bibitem {Tetsing}H.F. Tetsing, \textit{Normalized null hypersurfaces in
nonflat Lorentzian space forms satisfying }$L_{r}x=Ux+b$, Turk. J. Math.
\textbf{45} (2021), no. 4, 1809--1834.

\bibitem {Yang}B.G. Yang and X.M. Liu, \textit{Spacelike hypersurfaces in the
Lorentz-Minkowski space satisfying }$L_{r}x=Rx+b$\textit{,} Indian J. Pure
Appl. Math. \textbf{40} (2009), no. 6, 389-403.

\bibitem {Yang2}B. Yang and X. Liu,\textit{ Hypersurfaces satisfying }%
$L_{r}x=Rx$\textit{ in sphere }$S^{n+1}$\textit{ or hyperbolic space }%
$H^{n+1}$, Proc. Indian Acad. Sci. (Math. Sci.) \textbf{119} (2009), no. 4, 487--499.
\end{thebibliography}
\end{document}